\documentclass[reqno]{amsart}

\usepackage{amsmath, amsfonts, amsthm, amssymb, color,  graphicx, mathrsfs, cite}
\usepackage{comment,stmaryrd}

\theoremstyle{plain}
\newtheorem{theorem}{Theorem}[section]
\newtheorem{lemma}{Lemma}[section]
\newtheorem{proposition}{Proposition}[section]
\newtheorem{corollary}{Corollary}[section]

\theoremstyle{definition}
\newtheorem{definition}{Definition}[section]

\theoremstyle{remark}

\newtheorem{example}{Example}[section]

\numberwithin{equation}{section}
\allowdisplaybreaks

\usepackage{ifpdf}
\ifpdf \usepackage[colorlinks=true, citecolor=blue, linkcolor=blue, urlcolor=blue]{hyperref} \fi

\def\thrm{\begin{theorem}}
\def\thrml#1{\begin{theorem}\label{#1}}
\def\ethrm{\end{theorem}}
\def\lmm{\begin{lemma}}
\def\lmml#1{\begin{lemma}\label{#1}}
\def\elmm{\end{lemma}}
\def\dfntn{\begin{definition}}
\def\dfntnl#1{\begin{definition}\label{#1}}
\def\edfntn{\end{definition}}
\def\crllr{\begin{corollary}}
\def\crllrl#1{\begin{corollary}\label{#1}}
\def\ecrllr{\end{corollary}}
\def\xmpl{\begin{example}}
\def\xmpll#1{\begin{example}\label{#1}}
\def\exmpl{\end{example}}
\def\nmrt{\begin{enumerate}}
\def\enmrt{\end{enumerate}}
\def\qtn{\begin{equation}}
\def\qtnl#1{\begin{equation}\label{#1}}
\def\eqtn{\end{equation}}
\def\prpstn{\begin{proposition}}
\def\prpstnl#1{\begin{proposition}\label{#1}}
\def\eprpstn{\end{proposition}}

\def\tm#1{\item[{\rm (#1)}]}
\def\proof{{\bf Proof}.\ }
\def\eprf{\hfill$\square$}

\DeclareMathOperator{\aut}{Aut}

\DeclareMathOperator{\supp}{supp}

\newcommand{\Z}{\mathbb{Z}}

\newcommand{\ov}{\overline}

\def\qaq{\quad\text{and}\quad}

\def\lg{\langle}
\def\rg{\rangle}

\def\cZ{\mathcal {Z}}

\def\cA{\mathcal {A}}

\def\cK{{\mathcal K}}
\def\cP{{\mathcal P}}
\def\cD{{\mathcal D}}
\def\cB{{\mathcal B}}

\def\qaq{\quad\text{and}\quad}
\def\qoq{\quad\text{or}\quad}

\def\ov{\overline}

\begin{document}

\title{Schur Rings over $\cZ\times \cZ_3$}

\author{Gang Chen}
\address{School of Mathematics and Statistics, Central China Normal University, Wuhan, China}
\email{chengangmath@mail.ccnu.edu.cn}

\author{Jiawei He}
\address{School of Mathematics and Statistics, Central China Normal University, Wuhan, China}
\email{hjwywh@mails.ccnu.edu.cn}

\thanks{The authors are supported by the NSFC grant No. 11971189.}

\begin{abstract}
For the direct product $\cZ\times \cZ_3$ of infinite cyclic group $\cZ$ and a cyclic group $\cZ_3$ of order $3$, the schur rings over it are classified. In particular, all the schur rings are proved to be traditional. 

\end{abstract}


\date{}

\maketitle


\section{Introduction}\label{in}

Schur rings (or $S$-rings) were originally developed by Schur \cite{S33} and Wielandt \cite{W49} in order to study permuation groups as an alternative to character theory and have since been used in many other applications. Schur rings have been widely unsed in the study of association schemes, strongly regular graphs, and different sets.  For the applications of Schur rings to association schemes, see for example Sec. 2.4 and Sec. 4.4 of \cite{CP}.  \medskip

Classification of Shcur rings began with Schur himself, when Schur conjectured that all Schur rings over a group $H$ conincides with a so-called {\it transitivity module,} a submodule of the group ring $\Z[H]$ constructed by a group $G$ acting on its regular subgroup $H$. Such Schur rings are said to be {\it Schurian}. Wielandt \cite{W64} first showed that not all Schur rings are Schurian. Leung and Man in \cite{LM96} and \cite{LM98} showed that all Schur rings over a finite cyclic groups fall into one of four families: partitions induced by orbits of automorphism subgroups ({\it orbit Schur rings}), partitions induced by a group factorization ({\it tensor products}), partitions induced by cosets ({\it wedge products}), and the {\it trivial Schur ring}. We say that a Schur ring is {\it trational } if it is one of these four types.  Thus, the main results of Leung and Man tell us that all Schur rings over finite cyclic groups are traditional. Note that the above mentioned non-Schurian example by Wielandt is non-traditional; see the last paragraph of the fourth page of  \cite{M14}.  \medskip

Rcently, Bastian et al. \cite{B19}  study the Schur rings over infinite groups.  Schur rings over the integers, over a virtually infinite cyclic group, and over torsion-free locally cyclic groups are classified. All these Schur rings are traditional. However, there are non-triditional Schur rings over free groups and free products. \medskip 

Let $\cZ=\lg z\rg$ denote the infinite cyclic group, written multiplicative. The finite cyclic group of order $n$ will be denoted by $\cZ_n=\lg a\rg$. In the present paper, all Schur rings over the direct product $\cZ$ and $\cZ_3$ will be classified. The following is the main result. 

\thrml{1451b}All schur rings over the group $\mathcal{Z} \times \mathcal{Z}_{3}$ are of one of the following forms:
\nmrt
\tm{i}  $F[\cZ_3] \wedge F[\cZ]$ or $ F[\cZ_3] \wedge F[\cZ]^{\pm}$,

\tm{ii} $F[H \times \cZ_3] \wedge F[\cZ]$ or $F[H \times \cZ_3]^{\pm} \wedge F[\cZ]^{\pm}$,

\tm{iii} $F[K \times \cZ_3] \wedge F[\cZ]$ or $F[K \times \cZ_3]^{\pm} \wedge F[\cZ]^{\pm}$,

\tm{iv} $F[\cZ \times \cZ_3]$ or $F[\cZ \times \cZ_3]^{\pm}$,

\tm{v} an orbit Schur ring. 
\enmrt
where $H \leq  K\leq\cZ$.
Hence, all Schur rings over $\cZ \times \cZ_3$ are either orbit or wedge product Schur
rings, which implies that the group $\cZ\times \cZ_3$ is traditional.
\ethrm

Throughout let $G$ denote an arbitrary group with identity element $1$ and $F$ an arbitrary field of characteristic zero, which contains $\Z$ as a subring. The group algebra over $G$ with coefficients in $F$ is denoted $F[G]$. The torsion subgroup of $G$ is denoted $T(G)$.  
\section{Preliminary}\label{Pr}

 For any finite non-empty subset $C\subseteq G$, denote $\ov{C}=\sum_{g\in C}g$ and call $\ov{C}$ a {\it simple quantity}. In addition, we define $C^*=\{g^{-1}: g\in C\}$. A partition $\cP$ of $G$ is said to be  {\it finite support} provided any $C\in \cP$ is a finite subset of $G$. Our notation for Schur ring is taken from \cite{B19} and \cite{P15}. 
 
 \subsection{Schur Rings}
\dfntnl{1055a} Let $\cP$ be a finite support partition of $G$ and $\cA$ a subspace of $F[G]$ spanned by $\{\ov{C}: C\in \cP\}$. We say that $\cA$ is a {\it Schur ring} (or $S$-ring) over $G$ if  
 \nmrt
  \tm{1} $\{1\} \in \cP$, 
  \tm{2} for any $C\in \cP$, $C^*\in \cP$, 
  \tm{3} for all $C, D\in \cP$, $\ov{C}\cdot\ov{D}=\sum_{E\in \cP}\lambda_{CDE}\ov{E}$, where all but finitely many $\lambda_{CDE}$ equal $0$. 
 \enmrt
\edfntn
 
For an $S$-ring $\cA$ over $G$, the associated partition $\cP$ is denoted $\cD(\cA)$ and each element in $\cD(S)$ is called a {\it basic  set}
of $\cA$. We say that a subset $C$ of $G$ is an {\it $S$-set} if $C$ is a union of some basic sets of $\cA$. An {\it $S$-subgroup} is simutaneously an $S$-set and a subgroup.\medskip
 
Suppose $\alpha=\sum_{g\in G}\alpha_gg$ and $\beta=\sum_{g\in G}\beta_gg\in F[G]$ (note that here only finitely many nonzero coefficients $\alpha_g$ and $\beta_g$). Then define 
$\alpha^*=\sum_{g\in G}\alpha_gg^{-1}$ and $\supp(\alpha)=\{g: \alpha_g\ne 0 \}$. Additionally,  the {\it Hadamard product} is defined as $\alpha\circ \beta=\sum_{g\in G}\alpha_g\beta_gg$. \medskip

The following theorem was proved by Wielandt in \cite{W64} when $G$ is finite. The general case appeared as \cite[Corollary 2.12]{B19}. 
\thrml{1505b} Let $G$ be a group and $\cA$ a subring of $F[G]$. Then $\cA$ is a Schur ring if and only if $\cA$ is closed under $\circ$ and $*$, $1\in \cA$, and for all $g\in G$ there exists some $\alpha\in \cA$ such that $g\in \supp(\alpha)$. 
\ethrm

Let $f: F\rightarrow F$ be a function such that $f(0)=0$. For any $\alpha=\sum_{g\in G}\alpha_gg$, we set
$$
f(\alpha)=\sum_{g\in G}f(\alpha_g)g.
$$

The following proposition was proved by Wielandt in \cite[Proposition 22.3]{W64}. The general case appeared as \cite[Propostion 2.4]{B19}.
\prpstnl{913a}Let $\cA$ be a Schur ring over a group $G$ and the function be as above. Then 
$f(\alpha)\in \cA$ whenever $\alpha\in \cA$.
\eprpstn

Let $f$ take value $1$ at one non-zero number and $0$ at other numbers.  We get the following statement, which is known as the Schur-Wielandt principle; see \cite[Corollary 1.10]{P15}. 
 
\crllrl{1013c} Let $\cA$ be an $S$-ring over $G$. For any $\alpha\in \cA$, suppose $g\in \supp({\alpha})$ with $\alpha_g=c$. The set 
$$
\{g\in G: \alpha_g=c\}
$$
is an $S$-set. 
\ecrllr

The following two propositions, which were first proved by Wielandt in \cite{W64} and were generalized in \cite{B19},  tell us how to generate $S$-subgroups in any $S$-ring. 

\prpstnl{956b} Let $\cA$ be a $S$-ring over $G$. Let $\alpha\in \cA$ and 
$$
{\rm Stab}(\alpha)=\{g\in G: \alpha g=\alpha\}. 
$$
Then ${\rm Stab}(\alpha)$ is an $S$-subgroup of $\cA$.
\eprpstn

\prpstnl{956b} Let $\cA$ be an $S$-ring over $G$. Let $\alpha\in \cA$ and $H=\lg \supp(\alpha)\rg$. Then~$H$ is an $S$-subgroup of $\cA$. 

\eprpstn

Note that under the condition of Proposition \eqref{956b}, $\{C: C\in \cD(\cA), C\subseteq H\}$
consists of  basic sets of an $S$-ring over $H$. It is denoted $\cA_H$. 

When we have normal $S$-subgroup, we can construct $S$-ring over the factor group as shown in the following lemma. 

 \lmml{1440a}(\cite[Lemma 1.2]{LM90}) Let $\varphi: G\rightarrow H$ be a group homomorphism with ${\rm Ker}(\varphi)=K$ and $\cA$ be a Schur ring over $G$. Suppose that $K$ is an $S$-subgroup. Then the image $\varphi(\cA)$ is a Schur ring over $\varphi(G)$ where $\cD(\varphi(\cA))=\{\varphi(C): C\in \cD(S)\}$. 
\elmm

In particular, if $H=G/K$, the factor group of $G$ over $K$, and $\varphi$ is the natural homomorphism,  the corresponding $S$-ring is denoted $\cA_{G/K}$. 

\medskip

\subsection{Traditional Schur Rings}

The group ring $F[G]$ itself forms a $S$-ring over $G$ with basic set $\{ \{g\}: g\in G\}$. This is called the {\it discrete Schur ring} over $G$. When $G$ is a finite group, the partiotion $\{1, G\setminus 1\}$ produces a Schur ring known as the {\it trivial Schur ring} over $G$. Note that if $G$ is infinite there does not exist trivial $S$-ring over~$G$. \medskip

Assume that $G=H\times K$ and $\cA$ and $\cB$ are $S$-rings over $H$ and $K$, respectively. Then the set 
$$
\{CD: C\in \cD(\cA), D\in \cD(\cB)\}
$$
consists of a basic set  of an $S$-ring over $G$. The corresponding $S$-ring is called the {\it tensor  product} of $\cA$ and $\cB$. \medskip 

If $\cK$ is a finite subgroup of $\aut(G)$, the set of elements of $F[G]$ fixed by $\cK$ is an $S$-ring over $G$, denoted $F[G]^{\cK}$ and called the {\it orbit Schur ring} associated with $\cK$. Note that the discrete $S$-ring can be seen as an orbit $S$-ring by taking $\cK=1$.  \medskip 

An $S$-ring $\cA$ over $G$ is a {\it wedge product} if there exist nontrivial proper $S$-subgroups~$H,K$ such that $K\le H$, $K\unlhd G$, and every basic set outside  $H$ is a union of $K$-cosets. In this case, the series 
$$
1<K\le H<G
$$
is called a {\it wedge-decomposition} of $\cA$. Furthermore, we write $\cA=\cA_1\wedge \cA_2$, where 
$\cA_1=\cA_H$ and $\cA_2=\cA_{G/K}$. \medskip 

An $S$-ring over a group $G$ is called {\it traditional} if it is either a trivial Schur ring (when $G$ is finite), or a tensor product of Schur rings over smaller subgroups, and or a wedge product.  A group $G$ is called {\it traditional} if each $S$-ring over it is traditional. 

Several classes of groups are proved to be traditional as shown by the following theorem. 
\thrml{1556c} The following classes of groups are traditional:
\nmrt \tm{1}[\cite{LM96}, \cite{LM98}, \cite{B19}] cyclic groups. 
\tm{2}[\cite{B19}]$\cZ\times \cZ_2$. 
\tm{3}[\cite{B19}] Ttorsion-free locally cyclic groups. 
\enmrt
\ethrm

Actually, it is proved that the only Schur ring over $\cZ$ are discrete and symmetric Schur rings; see \cite[Theorem 3.3]{B19}. The symmetric Schur ring over $\cZ$ is denoted $F[\cZ]^{\pm}$.

\section{Schur rings over $\cZ\times \cZ_3$}

For $n\in \Z$, define the {\it $n$th Frobenius map} of $F[G]$ as 
	$$
	\alpha=\sum_{g\in G}\alpha_gg\, \mapsto \, \alpha^{(n)}=\sum_{g\in G}\alpha_gg^n.
	$$
 
 \lmml{1951b}(\cite[Theorem 23.9]{W64}, \cite[Theorem 2.20]{B19}) Suppose that $\cA$ is a Schur ring over $G$ where $G$ is an abelian group. Let  $m$ be an integer coprime to the orders of all torsion elements of $G$. Then for all $\alpha\in \cA$ we have $\alpha^{(m)}\in \cA$. In particular, if~$G$ is torsion-free, then all Schur rings over $G$ are closed under all Frobinius maps. 
 \elmm
 
 \lmml{1958c}(\cite[Lemma 2.21]{B19})Let $G$ be an abelian group such that the torsion subgroup $T(G)$ has finite exponent. Let $\cA$ be a Schur ring over $G$. Then $T(G)$ is an $S$-subgroup. 
 \elmm

The following result is known as the second Schur theorem on multipliers. 

\thrml{1535b} Let $\cA$ be an S-ring over an abelian group $G$, $p$ a prime divisor of~$|T(G)|$, and $E=\{g \in G: g^{p}=e\}$.  Then for every finite $S$ -set $X$ the set 
$$
X^{[p]}=\{x^{p}: x \in X,|X \cap E x| \not\equiv 0 (\bmod p)\}
$$ 
is an $S$-set.
\ethrm 
\proof
For a finite $S$-set $X$, one can see that

\qtnl{1042a}
\ov{X}^{p}=\left(\sum_{x \in X} x\right)^{p} \equiv \sum_{x \in X} x^{p}=\ov{{X}^{(p)}} \, (\bmod p),
\eqtn

where we write 
$$
\sum_{g} a_{g} g \equiv \sum_{g} a_{g}^{\prime} g\, (\bmod p)
$$
 if $a_{g} \equiv a_{g}^{\prime}(\bmod p)$ for all $g\in G$.  
 
 On the other hand, given a coset $ C \in G / E,$ the element $g^p$ does not depend on the choice of $g \in C$.  Denote it by $h_{C}$.  Then the mapping $C \mapsto h_{C}$ is a bijection from $G / E$ onto $G^{p}$.  So
\qtnl{1043b}
\ov{{X}^{(p)}}=\sum_{C \in G / E}|C \cap X| \cdot h_{C} \equiv \sum_{C \in G / E,|C \cap X| \not \equiv 0 \, (\bmod p)}|C \cap X| \cdot h_{C} \, (\bmod p).
\eqtn
Let $f: F \rightarrow F$ be the function such that $f(x)=1$ if $x$ is an integer not divisible by $p$; $f(x)=0$ for other $x\in F$. The formulas \eqref{1042a} and \eqref{1043b} imply that the set 
$$
\left\{h_{C}: C \in G / E,|C \cap X| \not \equiv 0(\bmod p)\right\}
$$ 
is an $S$-set. Since this set equals $X^{[p]}$,  we are done.
\eprf

\medskip

The following two propositions are useful in analyzing the structure of $S$-rings over $\cZ\times \cZ_p$, with $p$ an odd prime. 

\prpstnl{1501a}Let $\cA$ be an S-ring over $\cZ \times \cZ_p=\langle z\rangle\times \langle a\rangle$, and let $H$ be the unique maximal $\cA$-subgroup contained in $\cZ$. If there exists $z^ma^i \in X\in \cD(\cA)$ with~$|X|<p$, then $z^{pm}\in H$.
\eprpstn
\proof Using the notation in Theorem \ref{1535b}, we have $X^{[p]}=X^{(p)}$ by the assumption. Furthermore, it is an $S$-set and hence $\langle X^{(p)}\rangle$ is an $S$-subgroup. According to the maximality of $H$, one can see that
$$
(z^ma^i)^p\in X^{(p)}\subseteq H, 
$$
as required.
\eprf

\medskip

The following basic fact will be used freely in the proof of the main theorem. 

{\bf Remark.}\, {\it Let $\cA$ be an $S$-ring over $G$ and $K$ a normal  $S$-subgroup. Let $\varphi$ be the natural homomorphism from  $G$ onto $G/K$. For any $X\in \cD(\cA)$, $\varphi^{-1}(X)$ is an~$S$-set since it is equal to $XK$.}

\medskip

\thrml{1521b} Let $\cA$ be an S-ring over $\cZ \times \cZ_p=\lg z\rg\times \lg a\rg$, where $p$ is an odd prime. Then for all $X\in \cD(\cA)$ either $X=z^m\mathcal{Z}_{p}$ or $p\neq |X|$.

\ethrm
\proof
Denote $G=\cZ \times \cZ_p=\lg z\rg \times\lg a\rg$. Notice by Lemma \ref{1958c} that $T(G)=\cZ_p$ is an $S$-subgroup. 

According to Lemma~\ref{1440a}, $\cA_{G/T(G)}$ is a Schur ring over $G/T(G)$. 
Moreover,~$\cA_{G/T(G)}$ is either discrete or symmetric. 
By the above remark, it follows that 
$$
X\subseteq z^m\mathcal{Z}_{p} \:\:{\rm or}\:\: X\subseteq \{z^m,z^{-m}\}\mathcal{Z}_{p}, 
$$
for all $X\in \cD(\cA)$, where $z^m\in \cZ$ satisfies $X\cap z^m\cZ_p\neq \varnothing$. 
   
Towards a contradiction, assume that  there exists  $X\in \cD(\cA)$ with  $|X|=p$ and $X\neq z^m\cZ_{p}$ for all $m\in \Z$.  It follows that  $X\subseteq \{z^m,z^{-m}\}\cZ_{p}$ for some $m\in \Z$ such that 

\qtnl{1039b}
X\cap z^m\cZ_{p}\neq \varnothing \qaq X\cap z^{-m}\cZ_{p}\neq \varnothing.
\eqtn

So, we may assume that 
$$
X=\{z^ma^{i_1},z^ma^{i_2}, \ldots,z^ma^{i_t},z^{-m}a^{j_1},z^{-m}a^{j_2},\ldots,z^{-m}a^{j_k}\}. 
$$
 Here, ${i_1},\ldots,{i_t},{j_1},\ldots,{j_k}\in \mathbb{N}$ satisfying 
 
\qtnl{1612b} 
t+k=p \qaq  t\neq 0\neq k.  
\eqtn

 It follows that 
  $$\ov{X}(\ov{T(G)}-1)= t\overline{z^mT(G)}+k\ov{z^{-m}T(G)}-\ov{X}\in \cA.
$$
This implies that  $t\ov{z^mT(G)}+k\ov{z^{-m}T(G)}\in \cA$ and hence 
$$
t\ov{\{z^m,z^{-m}\}T(G)}+(k-t)\ov{z^{-m}T(G)}\in \cA. 
$$
By \eqref{1612b}, $k\neq t$ as $p$ is odd. Therefore,  $(k-t)\ov{z^{-m}T(G)}\in \cA$ yields that  
$\ov{z^{-m}T(G)}\in$~ $\cA$.  As a result,  $X\subseteq z^{-m}T(G)$. This is a contradiction  to \eqref{1039b}.  
\eprf

\medskip

{\bf Proof of Theorem \ref{1451b}}\, In the sequel, we fix the following notation: 
$$
G=\cZ \times \cZ_{3}=\lg z\rg \times\lg a\rg, 
$$
$\cA$ is  a Schur ring over $G$, and $H$ is  the unique maximal~$S$-subgroup contained in~$\cZ$.

\medskip 

Note that $T(G)=\{1, a,a^2\}=\cZ_{3}$ is an $S$-subgroup by Lemma \ref{1958c}. Let
 $$
 \pi: \cZ \times \cZ_{3} \rightarrow \cZ
 $$ 
 be the natural projection map, with  $\ker (\pi)=\cZ_{3}$.  Then $\pi(\cA)$ is a Schur ring over~$\cZ$,  which  is either discrete or symmetric by \cite[Theorem 3.3]{B19}.

 \medskip
 
{\bf Case 1.} \,  Suppose $\pi(\cA)=F[\cZ]$.  Then $\cA_H$ is discrete. Furthermore, 
 $$
 \left\{z^{m}, a z^{m}, a^2 z^{m}\right\}=\pi^{-1}\left(\left\{z^{m}\right\}\right)
 $$ 
 is an $S$-set for all integers $m$. Hence 
 $$
 \cA_{H\times \lg a\rg}=\cA_H\otimes \mathcal{A}_{ \lg a\rg}.
 $$
  If $H=\cZ$,  then $\cA=\cA_{\cZ}\otimes \cA_{ \lg a\rg}$ and hence the $S$-ring $\cA$ is traditional.
  
  \medskip
  
 Next, we assume that $H<\cZ$.     
 Choose  $z^{m} \in \cZ \setminus H$.   Then
 $$
 \{z^{m}, a z^{m}\}\in \cD(\cA), \, \, {\rm  or} \, \,  \{z^{m}, a^2 z^{m}\}\in \cD(\cA), \, \, {\rm  or} \, \, \{z^{m}, a z^{m}, a^2 z^{m} \}\in \cD(\cA).
 $$
Suppose 
$$
\{z^{m}, a z^{m}, a^2 z^{m}\}\in \cD(\mathcal{A}). 
$$  
In this case, we may assume that 
$$
\cA=\cA_{H\times \cZ_3} \wedge \cA_{G/\cZ_3}. 
$$ 
Otherwise, there exists $z^{m'} \in \cZ \setminus H$ such that
$$\{z^{m'}, a z^{m'}\}\in \cD(\cA), \qoq \{z^{m'}, a^2 z^{m'}\}\in \cD(\cA).
$$
We may replace $z^{m}$ with $z^{m'}$.
  
Since both $a$ and $a^2$ are generators of $\cZ_3$, without loss of generality we may assume that $\{z^{m}, a z^{m}\}\in \cD(\cA)$.   Therefore,  $a^2z^m\in \cD(\cA)$.  It follows that  $z^{3m}\in \cA$ by Proposition  \ref{1501a}. 
It is easy to see that $\{a, a^2\}\in \cD(\cA)$.   Let 
$$
K=\lg z^{m}, H\rg, 
$$
then $K$ is the unique subgroup of $G$ such that $|K: H|=3$.  Also, 
$$
K \times \cZ_{3}=\lg z^{m}, a z^{m}, H\rg 
$$
is an $S$-subgroup.

\medskip

{\bf Subcase 1.1.}\,  $K<\cZ$.  Then for any  $z^{n} \in \cZ \setminus K$,  the basic set containing $z^{n}$ must be $\{z^{n}, z^{n}a, a^2 z^{n}\}$. 
 Otherwise, we have $\{z^{3n}\}\in \cD(\cA)$ by Proposition \ref{1501a} and hence $z^{n} \in K$, a contradiction. 
 Thus, 
  $$
  \cA=\cA_{K\times \cZ_3}\wedge \cA_{G/\cZ_3}. 
  $$
 Consequently,  $\cA$ is  traditional whenever $1< K< \cZ$. 
 
 \medskip
   
{\bf Subcase 1.2.}\,  $K=\cZ$.  Then
$$
\{z, a z\}\in \cD(\cA), \, \, {\rm  or} \, \,  \{z, a^2 z\}\in \cD(\cA), \, \, {\rm  or} \, \, \{z, a z, a^2 z \}\in \cD(\cA). 
$$
Suppose 
\qtnl{1619c}
\{z, a z, a^2 z\}\in \cD(\cA). 
\eqtn 

\medskip

{\bf Claim 1.}\, {\it Under the assumption in \eqref{1619c}, we have $\cA=\cA_{H\times \cZ_3}\wedge \cA_{G/\cZ_3}$. In particular, $\cA$ is traditional.}

\proof If the claim is false, then there exists $3\nmid k\in \Z$ such that
$$
\{z^k, a z^k\}\in \cD(\cA),\qoq \{z^k, a^2 z^k\}\in \cD(\cA).
$$ 
    
If $k\equiv 1\:\:(mod\: 3)$, 
then
  $$
\{z,az\}=\{z^{1-k}\}\{z^k, a z^k\}\in \cD(\cA), \qoq \{z,a^2z\}=\{z^{1-k}\}\{z^k, a^2 z^k\}\in \cD(\cA), 
 $$ 
which  is a contradiction to statement \eqref{1619c}.  
      
 If $k\equiv 2\:\:(mod\: 3)$,  then
 $$
 \{z,a^2z\}=\{z^{1+k}\}\{z^{-k}, a^2 z^{-k}\}\in \cD(\cA) \:\: {\rm or} \:\:\{z,az\}=\{z^{1+k}\}\{z^{-k}, a z^{-k}\}\in \cD(\cA), 
 $$ 
a contradition.
 
It follows that $\cA=\cA_{H\times \cZ_3}\wedge \cA_{G/\cZ_3}$ and consequently $\cA$ is traditional.
 
\eprf
  
\medskip
        
 Now, we may assume that $\{z, a z\}\in \cD(\cA)$. One can see that 
 $$
 a^2z\in \cD(\cA) \qaq \lg z^3\rg\leq H\leq \cZ.
$$
 Observe that for any $k\in \Z^+$ and $3\nmid k$, 
 $$
  \{z^k,a^kz^k\}=\{a^{2k-2}z^k,a^{2k-1}z^k\}=(a^{2}z)^{k-1}\{z, a z\}\in \cD(\cA). 
 $$ 
Henc, the basic sets in $\cD(\cA)$  are the following: 
 $$
 \cD(\cA)=\{1\}; \{z^k,a^kz^k\}, \, k \in \Z,3\nmid k; \{a^{2k}z^k\}, \,  k \in \Z,3\nmid k; \{az^k,a^2z^k\}, \, k \in \Z,3| k.
$$
Let $\psi$ be the automorphism defined as follows: 
 \qtnl{1646c}
 \psi: G \rightarrow G,\:\:z\mapsto az,\:\:a\mapsto a^2.
 \eqtn
 It is easy to see that $\cA=F[G]^{\lg \psi\rg}$. Thus, $\cA$ is traditional. 
 
 \medskip  
 
 Similarly, if $\{z, a^2 z\}\in \cD(\cA)$,  then $\cA=F[G]^{\lg \delta \rg}$. Here,  $\delta \in \aut(G)$ defined as:
 \qtnl{1025b}
 \delta: G \rightarrow G,\:\:z\mapsto a^2z,\:\:a\mapsto a^2.
 \eqtn
 
 \medskip
 {\bf Case 2.}\, Suppose $\pi(\cA)=F[\cZ]^{\pm}$. Then, 
 $$
\{z^{m}, z^{-m}, a z^{m}, a z^{-m},a ^2z^{m}, a ^2 z^{-m}\}=\pi^{-1}\{z^{m}, z^{-m}\}
 $$
  is an $S$-set for all $m\in \Z$. 
  
  Observe that for any nonzero integer $m$, $\{a^kz^m, a^kz^{-m}\}$ can not be a basic set. Otherwise, $\{z^{3m}\}$ would be a basic set by Proposition \ref{1501a}.  This is a contradiction to the assumption.
  
  \medskip 
  
{\bf Subcase 2.1.}\, $H=\cZ$.  Thus, $\{z^m, z^{-m}\}\in \cD(\cA)$ for all $m\in \Z$. 
It follows that 
$$
\{a z^{m}, a z^{-m},a ^2z^{m}, a ^2 z^{-m}\}
$$
is an $S$-set for all $m\in \Z$. 

\medskip

If $\cA_{\lg a\rg}$ is discrete, then $\cA=\cA_{\cZ}\otimes\cA_{\lg a\rg}$ and so $\cA$ is traditional.

\medskip

Next, suppose that
$$
\{a,a^2\}\in \cD(\cA) \qaq \{az, az^{-1}, a^2z, a^2z^{-1}\}\in \cD(\cA).
$$ 
If there exists $k\in \mathbb{Z}^+$ such that $\{az^k,az^{-k}\}\in \cD(\cA)$, then 
$$
az^{2k+1}+az^{-1}+az+az^{-2k-1}=(az^k+az^{-k})(z^{k+1}+z^{-k-1})\in \cA.
$$
This implies that $\{az^{-1},az\}\in \cD(\cA)$,  a contradiction. Similarly, it is impossible that  $\{az^k,a^2z^{-k}\}\in \cD(\cA)$.   We conclude that 
$$
\{az^k,az^{-k},a^2z^k,a^2z^{-k}\}\in \cD(\cA), \forall  k\in \mathbb{Z}^+.  
$$
And therefore,  $\cA=\cA_{\cZ}\otimes\cA_{\lg a\rg}$.  In particular, $\cA$ is traditional. 

\medskip

Now, suppose that 
$$
\{a,a^2\}\in \cD(\cA) \qaq \{az,a^2z^{-1}\}
\in \cD(\cA).
$$
By the observation at the beginning of Case 2, $\{az^{-1}, a^2z\}$ is also a basic set. 

Applying Lemma \ref{1951b} to $\{az,a^2z^{-1}\}$ and $\{az^{-1}, a^2z\}$, we obtain that 
$$
\{a^kz^{-k}, a^{2k}z^k\}\qaq \{a^kz^k, a^{2k}z^{-k}\}\in \cD(\cA), \forall  3\nmid k\in \Z.  
$$
Thus, for any $3\nmid k$,  
$$
a^{k+1}z^{k-1}+a^{2k+1}z^{-k-1}+a^{k+2}z^{k+1}+a^{2(k+1)}z^{-k-1}=(az^{-1}+a^2z)(a^kz^k+a^{2k}z^{-k})\in \cA. 
$$
In particular, if $k\equiv 1\:\:(mod\: 3)$, one can see that 
$$
\{az^k, a^2z^{-k}\}\in \cD(\cA), \forall 3|k. 
$$
We conclude that $\{az^k, a^{2}z^{-k}\}\in \cD(\cA)$ for all $ k\in \Z$.  It is easy to see that 
$$
\cA=F[G]^{\lg\xi\rg},
$$
where $\xi\in \aut(G)$ is defined as
\qtnl{1649b}
\xi:  G \rightarrow G, \:\:z\mapsto z^{-1},\:\: a\mapsto a^2.
\eqtn

\medskip

Finally, neither  $\{az, az^{-1}\}$ nor  $\{az, a^2z\}$ is a basic set by the observation at the beginning of Case 2. 

\medskip 

{\bf Subcase 2.2.}\, $H< \cZ$. 

Choose $z^{m} \in \cZ \setminus H$. Then one of the following holds:
\nmrt 
\tm{i} $\{z^{m}, z^{-m}, a z^{m}, a z^{-m},a ^2z^{m}, a ^2 z^{-m}\}\in \cD(\cA)$,
\tm{ii} $\{z^{m}, z^{-m},a z^{m},a ^2 z^{-m}\}\in \cD(\cA)$ and  $\{a z^{-m}, a^2 z^{m}\}\in \cD(\cA)$, 
\tm{iii}$\{z^{m}, z^{-m},a z^{-m},a ^2 z^{m}\}\in \cD(\cA)$ and $\{a^2 z^{-m}, a z^{m}\}\in \cD(\cA)$, 
\tm{iv}$\{z^{m},  a z^{m}, a z^{-m}\}\in \cD(\cA)$ and $\{ z^{-m}, a^2 z^{-m}, a ^2z^{m}\}\in \cD(\cA)$,
\tm{v} )$\{z^{m},  a^2 z^{-m}, a z^{-m}\}\in \cD(\cA)$ and $\{ z^{-m}, a z^{m}, a ^2z^{m}\}\in \cD(\cA)$, 
\tm{vi}$\{z^{m},  a^2 z^{-m}, a^2 z^{m}\}\in \cD(\cA)$ and $\{ z^{-m}, a z^{-m}, a ^2z^{m}\}\in \cD(\cA)$,
\tm{vii}$\{z^{m}, a z^{-m}\}\in \cD(\cA)$, $\{ z^{-m}, a ^2z^{m}\}\in \cD(\cA)$ and $\{ a^2z^{-m}, a z^{m}\}\in \cD(\cA)$, 
\tm{viii}$\{z^{m}, a^2 z^{-m}\}\in \cD(\cA)$, $\{ z^{-m}, az^{m}\}\in \cD(\cA)$ and $\{ az^{-m}, a^2 z^{m}\}\in \cD(\cA)$.
\enmrt

If (i)  holds for all $z^m\in \cZ\setminus H$ , then 
$$
\mathcal{A}=\mathcal{A}_{H\times \mathcal{Z}_3}\wedge\mathcal{A}_{G/\mathcal{Z}_3}.
$$
In particular, $\cA$ is traditional. 

\medskip 

Conditions (iv),  (v),  and (vi) can not happen  by Theorem \ref{1521b}. 

\medskip

 Now, assume that statement (ii), (iii), (vii),   or (viii) holds. 
In these cases, we conclude that $z^{3m}\in H$ by Proposition \ref{1501a}. Set
 $$
 M=\lg z^{m}, H\rg=\lg z^{-m}, H\rg.
 $$ 
 Then $M$ is the unique subgroup of $G$ such that $|M: H|=3$. 
 In addition, 
 $$
 M \times \cZ_{3}=\lg z^{m},z^{-m}, a z^{m},a^2z^{-m}, H\rg=\lg z^{m},z^{-m}, a z^{-m},a^2z^{m}, H\rg
 $$
is an $S$-subgroup. 
  
If $M<\mathcal{Z}$, then for any $z^{l} \in \mathcal{Z} \backslash M$, by the above arguments, one can easily see that  the basic  set containing $z^{l}$ must be 
$$
\{z^{l}, z^{-l}, a z^{l}, a z^{-l},a ^2z^{l}, a ^2 z^{-l}\}.
   $$
  Hence,  
  $$
  \cA=\cA_{M\times \cZ_3}\wedge \cA_{G/\cZ_3}, 
  $$
 and so $\cA$ is traditional.
 
 \medskip
   
If $M=\cZ$,  then $H=\lg z^3\rg$ by $|M:H|=3$. There are fives cases to consider. 
\nmrt
\tm{1}$\{z, z^{-1}, a z, a z^{-1},a ^2z, a ^2 z^{-1}\}\in \cD(\cA)$, 
\tm{2}$\{z, z^{-1},a z,a ^2 z^{-1}\}\in \cD(\cA)$ and $\{a z^{-1}, a^2 z\}\in \cD(\cA)$, 
\tm{3}$\left\{z, z^{-1},a z^{-1},a ^2 z\right\}\in \cD(\cA)$,  and $\{a^2 z^{-1}, a z\}\in \cD(\cA)$, 
\tm{4}$\{z, a z^{-1}\}\in \cD(\cA)$, $\{ z^{-1}, a ^2z\}\in \cD(\cA)$,  and $\{ a^2z^{-1}, a z\}\in \cD(\cA)$,
\tm{5}$\{z, a^2 z^{-1}\}\in \cD(\cA)$, $\{ z^{-1}, az\}\in \cD(\cA)$,  and $\{ az^{-1}, a^2 z\}\in \cD(\cA)$.
\enmrt

\medskip

Suppose condition (1) holds. Then we have the following claim. 

\medskip

{\bf Claim 2.}\,  $\cA=\cA_{H\times \cZ_3}\wedge \cA_{\cZ_3}$.

\proof If the claim is false, then there exists $k\in \Z^+$ not divisible by $3$ such that 
$$
\{a^2z^{-k}, az^k\}\in \cD(\cA) \qoq  \{az^{-k}, a^2z^k\}\in \cD(\cA). 
$$
First, assume that $\{a^2z^{-k}, az^k\}\in \cD(\cA)$. 

\medskip

If $k\equiv 1\:\:(mod\: 3)$,  then $\{z^{k-1},z^{1-k}\}\in \cD(\cA)$ by $3|{k-1}$. It follows that 
$$
a^2z^{-1}+ a^2z^{1-2k}+ az^{2k-1}+ az=(a^2z^{-k}+ az^k)\cdot (z^{k-1}+z^{1-k})\in \cA.
$$ 
This yields  that $\{a^2z^{-1},az\}\in \cD(\cA)$, a contradiction. 

\medskip

If $k\equiv 2\:\:(mod\: 3)$,  then $\{z^{k+1},z^{-1-k}\}\in \cD(\cA)$ by $3|{k+1}$.  It follows that  
$$
az^{-1}+a^2z^{-1-2k}+ az^{2k+1}+ a^2z=(a^2z^{-k}+ az^k)\cdot (z^{k+1},z^{-1-k})\in \cA.
$$ 
Thus,  $\{az^{-1},a^2z\}\in \cD(\cA)$, a contradiction. 

\medskip

Therefore,  the claim is valid in this case. Similarly, the claim holds if~$\{az^{-k}, a^2z^k\}$ is a basic set 
for some positive integer  $k$ not divisible by $3$. 
\eprf

\medskip

Suppose condition (2) holds. Then, by Lemma \ref{1951b} $\{a^2z^{-2},az^2\}=\{az^{-1},a^2z\}^{(2)}$ is a basic set.  By the statements at the beginning of Subcase 2.2, this implies that 
$$
\{z^2, z^{-2},a z^{-2},a ^2 z^2\}\in \cD(\cA) \qoq \{z^2, a z^{-2}\}\in \cD(\cA).
$$
If $\{z^2, a z^{-2}\}\in \cD(\cA)$, then 
$$
az+a^2z^{-3}+a^2z^3+z^{-1}=(az^{-1}+a^2z)\cdot(z^2+az^{-2})\in \cA.
$$
 It means that $\{az,z^{-1}\}$ is an $S$-set,  a contradiction to the assumption of (2). We conclude that 
 \qtnl{1455a}
 \{z^2, z^{-2},a z^{-2},a ^2 z^2\}\in \cD(\cA).
 \eqtn
 Note that 
 $$
 \{az^3, az^{-3}\}\in \cD(\cA), \qoq \{az^3, a^2z^{-3}\} \in \cD(\cA), \qoq \{az^3, az^{-3}, a^2z^{-3}, a^2z^3\}\in \cD(\cA). 
 $$
 
 If $\{az^3, az^{-3}\}\in \cD(\cA)$,  then 
 $$
 a^2z^2+z^4+a^2z^{-4}+z^{-2}=(az^3+ az^{-3})(az^{-1}+a^2z)\in \cA.
 $$
So, $\{a^2z^2,z^{-2}\}$ is an $S$-set, a contradiction to statement \eqref{1455a}.

 If $\{az^3, a^2z^{-3}\} \in \cD(\cA)$,  then 
 $$
 a^2z^2+z^4+z^{-4}+az^{-2}=(az^3+ a^2z^{-3})(az^{-1}+a^2z)\in \cA
 $$
 and hence $\{a^2z^2,  az^{-2}\}$ is an $S$-set,  another  contradiction to statement \eqref{1455a}.
 
As a consequence, we obtain
\qtnl{1509b}
\{az^3, az^{-3}, a^2z^{-3}, a^2z^3\}\in \cD(\cA).
\eqtn

Next, we will prove the following claim. 

\medskip

{\bf Claim 3.}\, {\it The following statements hold:
\nmrt
\tm{a} $\{z^n, z^{-n}, a^nz^{n}, a^{2n}z^{-n}\}\in \cD(\cA)$, if  $n\in \Z^+$ and $3\nmid n$, 
\tm{b} $\{az^n, az^{-n}, a^2z^{n}, a^2z^{-n}\}\in \cD(\cA)$,  if  $n\in \Z^+$  and  $3| n$. 
\enmrt
}
\proof 
We prove the claim by induction on $n$. If $n=1$, $n=2$, or $n=3$, the claims follow by the assumption and statements \eqref{1455a} and \eqref{1509b}.

\medskip 

If $n\equiv 0\:\:(mod\: 3)$,  then $n-1\equiv 2\:\:(mod\: 3)$.  Hence, 
$$
\{z^{n-1},z^{1-n},az^{1-n},a^2z^{n-1}\}\in \cD(\cA),
$$
by inductive hypotheses.  

Assume that $\{az^n,az^{-n}\}\in \cD(\cA)$. Then 
$$
a^2z^{n-1}+z^{n+1}+a^2z^{-n-1}+z^{1-n}=(az^n+az^{-n})(az^{-1}+a^2z)\in \cA.
$$
It follows that $\{a^2z^{n-1},z^{1-n}\}$ is an $S$-set,  a contradiction. 

Assume that $\{az^n,a^2z^{-n}\}\in \cD(\cA)$. Then 
$$
a^2z^{n-1}+z^{n+1}+z^{-n-1}+az^{1-n}=(az^n+a^2z^{-n})\cdot (az^{-1}+a^2z)\in \cA.
$$
Thus,  $\{a^2z^{n-1},az^{1-n}\}$ is an $S$-set,  a contradiction. 
Hence, 
$$
\{az^n, az^{-n}, a^2z^{n}, a^2z^{-n}\}\in \cD(\cA). 
$$

\medskip 

If $ n\equiv 1\:\:(mod\: 3)$,  then $n-1\equiv 0\:\:(mod\: 3)$.  So, 
$$
\{az^{n-1}, az^{1-n},a^2z^{1-n}, a^2z^{n-1}\}\in \cD(\cA). 
$$
by inductive hypotheses. By Lemma \ref{1951b}, one can see that  
$$
\{az^{-n}, a^2z^n\}=\{az^{-1}, a^2z\}^{(n)}
$$
is an $S$-set.  Therefore,  $\{az^{-n}, a^2z^n\}\in \cD(\cA)$.   It  yields that 
$$
\{z^{n}, z^{-n},a z^{n},a ^2 z^{-n}\}\in \cD(\cA), \qoq  \{z^{n}, a^2 z^{-n}\}\in  \cD(\cA).
$$
Assume that $\{z^{n}, a^2 z^{-n}\}\in \cD(\cA)$.  Then 
$$
az^{n-1}+a^2z^{n+1}+z^{-n-1}+az^{1-n}=(z^{n}+a^2 z^{-n})\cdot(az^{-1}+a^2z)\in \cA.
$$
This implies that  $\{az^{n-1},az^{1-n}\}$ is an $S$-set,  a contradiction.  

Thus, 
$$
\{z^{n}, z^{-n}, a z^{n}, a ^2 z^{-n}\}\in \cD(\cA). 
$$

\medskip

If $n\equiv 2\:\:(mod\: 3)$,  then $n-1\equiv 1\:\:(mod\: 3)$.  Hence, 
 $$
 \{z^{n-1},z^{1-n},a^2z^{1-n},az^{n-1}\}\in \cD(\cA)
 $$ 
 by inductive hypotheses. 
 
 Similarly, we have $\{a^2z^{-n}, az^n\}\in \cD(\cA)$. One can see that 
 $$
 \{z^{n}, z^{-n},a z^{-n}, a ^2 z^{n}\}\in \cD(\cA), \qoq \{z^{n}, a z^{-n}\}\in \cD(\cA).
 $$ 
 Suppose that $\{z^{n}, a z^{-n}\}\in \cD(\cA)$,  then 
 $$
 az^{n-1}+a^2z^{n+1}+a^2z^{-n-1}+z^{1-n}=(z^{n}+ a z^{-n})\cdot(az^{-1}+a^2z)\in \cA.
 $$
This yields that $\{az^{n-1},z^{1-n}\}\in \cD(\cA)$,  a contradiction.  Hence, 
$$
\{z^{n}, z^{-n}, az^{-n}, a ^2 z^{n}\}\in \cD(\cA).
$$
This completes the proof of Claim 3.
\eprf

\medskip

Obviously, $\{a,a^2\}\in S(\mathcal{A})$.  As a consequence of Claim 3, one can easily see that 
$$
\cA=F[G]^{\lg \psi, \xi\rg}, 
$$ 
where $\psi$ and $\xi$ are automorphisms of $G$ defined as in \eqref{1646c} and \eqref{1649b}, respectivelty.

\medskip 

Suppose condition (3) holds.  One can similarly get that
$$
\cA=F[G]^{\lg \delta, \xi\rg},
$$
where $\delta$ and $\xi$ are defined as in \eqref{1025b} and \eqref{1649b}, respectively. 

\medskip

Suppose condition (4) holds. Applying Lemma \ref{1951b},  for all integers $k$ with $ 3\nmid k$  we obtain that
$$
\{z^k, a^k z^{-k}\}\in \cD(\cA),  \quad \{ z^{-k}, a ^{2k}z^k\}\in \cD(\cA),  \qaq \{ a^{2k}z^{-k}, a^k z^k\}\in \cD(\cA). 
$$
Assume that  $k\equiv 1\:\:(mod\: 3)$.  Observe that 
$$
az^{k-2}+a^2z^{k+2}+a^2z^{-k-2}+z^{2-k}=(z^k+a z^{-k})(az^{-2}+a^2z^{2})\in \cA. 
$$
This yields that  $\{a^2z^{k+2},a^2z^{-k-2}\}\in \cD(\cA)$.

Assume that  $k\equiv 2\:\:(mod\: 3)$.  Then 
$$
az^{k+1}+a^2z^{k-1}+z^{1-k}+az^{-1-k}=(z^k+a^2 z^{-k})(az+a^2z^{-1})\in \cA. 
$$
This implies that $\{az^{k+1},az^{-k-1}\}\in \cD(\cA)$. 

We conclude that 
$$
\{az^{k},az^{-k}\}\in \cD(\cA), \qaq \{a^2z^{k},a^2z^{-k}\}\in \cD(\cA), \forall  3| k\in \Z^+.
$$
Note that 
$$
(z+az^{-1})(z+az^{-1})=z^2+2a+az^{-2}\in \cA.
$$
Thus, $\{a\}\in \cD(\cA)$, and hence  $\mathcal{A}_{\langle a\rangle}$ is discrete.

Now let $\rho\in \aut(G)$ be defined as follows: 
\qtnl{1105c}
\rho: G \rightarrow G,\:\:z\mapsto az^{-1},\:\: a\mapsto a.
\eqtn

It is easy to see that 
$$
\cA=F[G]^{\lg \rho \rg}.
$$
In particular, $\cA$ is traditional. 

\medskip 

Suppose condition (5) holds.  Let $\sigma\in \aut(G)$ be defined as follows:
\qtnl{1110a}
\sigma: G \rightarrow G, \:\:z\mapsto a^2z^{-1},\:\:a\mapsto a.
\eqtn
Then one can see that 
$$
\cA=F[G]^{\lg \sigma \rg}.
$$

\end{document}